\documentclass[12pt]{article}
\usepackage{amsfonts,amssymb,latexsym}

\newtheorem{theorem}{Theorem}[section]
\newtheorem{lemma}[theorem]{Lemma}

\newtheorem{proposition}[theorem]{Proposition}

\newenvironment{proof}
{\par\addvspace{0.3cm}\noindent{\rm Proof. }}
{\nopagebreak\mbox{}\hfill $\Box$\par\addvspace{0.25cm}}
\newenvironment{proofof}[1]
{\par\addvspace{0.3cm}\noindent{\rm Proof of #1. }}
{\nopagebreak\mbox{}\hfill $\Box$\par\addvspace{0.25cm}}
\newcommand{\be}{\begin{equation}}
\newcommand{\ee}{\end{equation}}

\newcommand{\bq}{\begin{eqnarray}}
\newcommand{\eq}{\end{eqnarray}}
\newcommand{\nn}{\nonumber}
\newcommand{\ba}{\begin{array}}
\newcommand{\ea}{\end{array}}

\renewcommand{\Re}{\mathrm{Re\,}}

\newcommand{\even}{\mathrm{even}}
\newcommand{\wind}{\mbox{\rm wind\,}}

\newcommand{\wt}{\widetilde}

\newcommand{\iv}{^{-1}}
\newcommand{\iy}{\infty}
\newcommand{\ovl}{\overline}

\newcommand{\R}{{\mathbb R}}
\newcommand{\C}{{\mathbb C}}
\newcommand{\Z}{{\mathbb Z}}
\newcommand{\T}{{\mathbb T}}

\newcommand{\cR}{\mathcal{R}}

\newcommand{\om}{\omega}
\newcommand{\eps}{{\varepsilon}}
\renewcommand{\rho}{{\varrho}}
\renewcommand{\kappa}{\varkappa}

\newcommand{\Hp}{H^p(\T)}
\newcommand{\Hq}{H^q(\T)}
\newcommand{\Lp}{L^p(\T)}
\newcommand{\Lq}{L^q(\T)}
\newcommand{\Lpe}{L^p_{\rm even}(\T)}
\newcommand{\Lqe}{L^q_{\rm even}(\T)}
\newcommand{\Li}{L^\iy(\T)}

\newcommand{\Ci}{C^{\iy}(\T)}
\newcommand{\Cio}{C^{\iy}_{0}(\T)}

\begin{document}

\date{}
\title{Factorization of a class of Toeplitz + Hankel operators and
the $A_p$-condition}
\author{Estelle L. Basor\thanks{ebasor@calpoly.edu.
           Supported in part by NSF Grant DMS-0200167.}\\
                Department of Mathematics\\
                California Polytechnic State University\\
                San Luis Obispo, CA 93407, USA
         \and
         Torsten Ehrhardt\thanks{tehrhard@mathematik.tu-chemnitz.de.}\\
        		Fakult\"{a}t f\"{u}r Mathematik\\
          	Technische Universit\"{a}t Chemnitz\\
          	09107 Chemnitz, Germany}
\maketitle

\begin{abstract}
Let $M(\phi)=T(\phi)+H(\phi)$ be the Toeplitz plus Hankel operator
acting on $H^p(\T)$ with generating function $\phi\in L^\iy(\T)$.
In a previous paper we proved that $M(\phi)$ is invertible if and
only if $\phi$ admits a factorization $\phi(t)=\phi_{-}(t)\phi_{0}(t)$  
such
that $\phi_{-}$ and $\phi_{0}$ and their inverses belong to certain
function spaces and such that a further condition formulated in terms of
$\phi_{-}$ and $\phi_{0}$ is satisfied. In this paper we  prove that  
this additional condition is equivalent to the
Hunt-Muckenhoupt-Wheeden condition (or, $A_{p}$-condition) for a certain
function $\sigma$ defined on $[-1,1]$, which is given in terms of
$\phi_{0}$. As an application, a necessary and sufficient criteria for
the invertibility of $M(\phi)$ with piecewise continuous functions  
$\phi$ is
proved directly. Fredholm criteria are obtained as well.
\end{abstract}


\section{Introduction}
This paper is devoted to continuing the study
(started in \cite{BaEh2}) of operators of the form
\bq
  M(\phi) &=& T(\phi)+H(\phi)
  \eq
acting on the Hardy space $\Hp$ where $1<p<\iy$.
Here $\phi\in\Li$ is a Lebesgue measurable and essentially
bounded function on the unit circle $\T$.
The Toeplitz and Hankel operators are defined by
\be\label{f.2}
T(\phi):f\mapsto P(\phi f),\qquad
H(\phi):f\mapsto P(\phi (Jf)),
\qquad f\in\Hp,
\ee
where is $J$ the following flip operator,
\be\label{f.3}
J:f(t)\mapsto t\iv f(t\iv),\quad t\in\T,
\ee
acting on the Lebesgue space $\Lp$.
The operator $P$ stands for the Riesz projection,
\be
P:\sum_{n=-\iy}^\iy f_nt^n\mapsto \sum_{n=0}^\iy f_nt^n, \quad t\in\T,
\ee
which is bounded on $L^p(\T)$, $1<p<\iy$, and whose image is $H^p(\T)$.
The complex conjugate Hardy space $\ovl{\Hp}$ is the set of all
functions $f$ whose complex conjugate belongs to $\Hp$.
Moreover, we denote by $\Lpe$ the subspace of $\Lp$ consisting of all
even functions, i.e., functions $f$ for which $f(t)=f(t\iv)$.

For $\phi\in\Li$ we denote by $L(\phi)$ the multiplication operator
acting on $\Lp$,
\bq
L(\phi)&:&f(t)\mapsto \phi(t)f(t).
\eq
Obviously, $T(\phi)$ and $H(\phi)$ can be written as
$$
T(\phi)\;=\; PL(\phi)P|_{\Hp},\qquad
H(\phi)\;=\; PL(\phi)JP|_{\Hp}.
$$

In a previous paper \cite{BaEh2}  we proved that for
$\phi\in L^{\iy}(\T)$ the operator
$M(\phi)$ is a Fredholm operator on the space $H^{p}(\T)$ if and only
if the functions $\phi$ admits a certain kind of generalized
factorization. Before recalling the underlying definitions,
let us state the following simple necessary condition
for the Fredholmness of $M(\phi)$ which was also established in
\cite[Prop.~2.2]{BaEh2}. Therein $G\Li$ stands for the group of all  
invertible
elements in the Banach algebra $\Li$.

\begin{proposition}\label{p1.1}
Let $1<p<\iy$ and $\phi\in \Li$. If $M(\phi)$ is Fredholm on
$\Hp$, then $\phi\in G\Li$.
\end{proposition}

A function $\phi\in L^{\iy}(\T)$ is said to admit a {\em weak
asymmetric factorization in $\Lp$} if it can be written in the form
\bq\label{f:fact}
\phi(t) &=& \phi_{-}(t) t^{\kappa}\phi_{0}(t), \qquad t\in\T,
\eq
such that $\kappa\in\Z$ and
\begin{itemize}
\item[(i)]
$(1+t\iv)\phi_{-}\in\ovl{\Hp}$, $(1-t\iv)\phi_{-}\iv\in\ovl{\Hq}$,
\item[(ii)]
$|1-t|\phi_{0}\in\Lqe$, $|1+t|\phi_{0}\iv\in\Lpe$.
\end{itemize}
Here $1/p+1/q=1$. It was proved in \cite[Prop.~3.1] {BaEh2} that if a  
factorization
exists, then the {\em index} $\kappa$
of the weak asymmetric factorization is uniquely determined and the  
factors
$\phi_{-}$ and $\phi_{0}$ are uniquely determined up to a  
multiplicative constant.
(In \cite{BaEh2} also the notion of a {\em weak antisymmetric  
factorization in $\Lp$}
has been introduced. This notion will play no role in the present  
paper.)

In order to introduce yet another notion, let $\cR$ stand for the set  
of all
trigonometric polynomials. Under the assumption that $\phi$ admits a  
weak
asymmetric factorization in $\Lp$ introduce the linear spaces
\bq
\label{f.X1}
X_{1} &=& \Big\{\;(1-t\iv)f(t)\;:\;f\in\cR\;\Big\},\\[.5ex]
X_{2} &=& \Big\{\;(1+t\iv)\phi_{0}\iv(t)f(t)\;:\;f\in\cR,\;  
f(t)=f(t\iv)\;\Big\}.
\label{f.X2}
\eq
It is easy to see that $X_{1}$ and $X_{2}$ are linear subspaces of
$\Lp$ and that the space $X_{1}$ is dense in $\Lp$. Moreover,
it has been proved \cite[Lemma 4.1(a)]{BaEh2} that
\bq
B&:=& L(\phi_{0}\iv) (I+J)PL(\phi_{-}\iv)
\eq
is a well-defined linear (not necessarily bounded) operator acting from
$X_{1}$ into $X_{2}$.

We will call the above factorization (\ref{f:fact}) of $\phi$ an
{\em asymmetric factorization in $\Lp$} if in addition to (i) and (ii)  
the
following condition is satisfied:
\begin{itemize}
\item[(iii)]
The operator $B=L(\phi_{0}\iv) (I+J)PL(\phi_{-}\iv)$ acting from $X_{1}$
into $X_{2}$ can be extended by continuity to a linear bounded operator
acting on $\Lp$.
\end{itemize}
Clearly, due to the density of $X_{1}$ in $\Lp$ an equivalent  
formulation
for condition (iii) is the following statement:
\begin{itemize}
\item[(iii*)]
There exists a constant $M$ such that
$\|Bf\|_{\Lp}\le M\|f\|_{\Lp}$ for all $f\in X_{1}$.
\end{itemize}

The main result proved in \cite[Thm.~6.4]{BaEh1} is the following:
\begin{theorem}\label{t1.1}
Let $1<p<\iy$ and $\phi\in G\Li$. The operator $M(\phi)$ is a Fredholm  
operator on
$\Hp$ if and only if the function $\phi$ admits an asymmetric
factorization in $\Lp$. In this case, the defect numbers are given by
\bq\label{ker-coker}
\dim\ker M(\phi)=\max\{0,-\kappa\},\qquad
\dim\ker M^*(\phi)=\max\{0,\kappa\},
\eq
where $\kappa$ is the index of the factorization of $\phi$.
\end{theorem}

To formulate the main result of this paper we need the
notion of the Hunt-Muckenhoupt-Wheeden condition
(or, $A_p$-condition) with respect to the interval $[-1,1]$.

Let $1<p<\iy$, $1/p+1/q=1$, and let $\sigma:[-1,1]\to \R_+$
be a Lebesgue measurable, almost everywhere nonzero function.
Assume in addition that $\sigma\in L^{p}[-1,1]$ and
$\sigma\iv \in L^{q}[-1,1]$. We say that $\sigma$ satisfies the
{\em $A_{p}$-condition on $[-1,1]$} if
\bq
\sup_{I}\frac{1}{|I|}
\left(\int_{I}\sigma^{p}(x)\,dx\right)^{1/p}
\left(\int_{I}\sigma^{-q}(x)\,dx\right)^{1/q}
&<&\iy,
\eq
where the supremum is taken over all subintervalls $I$ of $[-1,1]$.
The length of the interval $I$ is denoted by $|I|$.
There is an intimate connection between the $A_{p}$-condition and  
the boundedness of the singular integral operator, which will be
stated later on.

The main result of this paper is the following:
\begin{theorem}\label{t1.3}
Let $1<p<\iy$, $1/p+1/q=1$, and $\phi\in \Li$.
The operator $M(\phi)$ is a Fredholm operator on $\Hp$ if and only if
the following conditions are satisfied:
\begin{itemize}
\item[(a)]
$\phi\in G\Li$.
\item[(b)]
The function $\phi$ admits a weak asymmetric factorization in
$\Lp$,
\bq
\phi(t)&=&\phi_{-}(t)t^{\kappa}\phi_{0}(t),\quad t\in\T.
\eq
\item[(c)]
The function
\bq\label{f1.sigma}
\sigma(\cos\theta) &:=&
|\phi_0\iv(e^{i\theta})|\,
\frac{|1+\cos\theta|^{1/2q}}{|1-\cos\theta|^{1/2p}}
\eq
satisfies  the $A_{p}$-condition.
\end{itemize}
Moreover, formulas (\ref{ker-coker}) hold in this case.
\end{theorem}
We note that it is straightforward to prove that condition (b) of the  
previous theorem implies
$\sigma\in L^p[-1,1]$ and $\sigma\iv\in L^q[-1,1]$.


\section{Proof of Theorem \ref{t1.3}}

Before we are able to give the proof of Theorem \ref{t1.3} we establish
some definitions and auxiliary results.

Let $C^{\iy}[-1,1]$ stand for the set of all infinitely differentiable
functions $f:[-1,1]\to\C$, and denote by $C^{\iy}_{0}[-1,1]$ the
subspace of all functions $f\in C^{\iy}[-1,1]$ such that
$f(x)$ and all of its derivatives vanish at the endpoints $x=-1$ and
$x=1$,
\bq
C^{\iy}_{0}[-1,1] &=&
\Big\{\;f\in C^{\iy}[-1,1]\;:\;f^{(n)}(-1)=f^{(n)}(1)=0
\mbox{ for all } n\ge0\;\Big\}.
\eq
The singular integral operator $S_{[-1,1]}$ is defined by the rule
\bq
(S_{[-1,1]}f)(x) &=& \frac{1}{\pi i}\int_{-1}^{1}\frac{f(y)}{y-x}\,dy,
\qquad x\in[-1,1],\
\eq
where the integral has to be inderstood as the Cauchy principle value.
For $f\in C^{\iy}_{0}[-1,1]$ the integral exists for each $x\in[-1,1]$.
In fact,
\bq
(S_{[-1,1]}f)(x) &=& \frac{1}{\pi  
i}\int_{-1}^{1}\frac{f(y)-f(x)}{y-x}\,dy+
\frac{f(x)}{\pi i}\ln\left(\frac{1-x}{1+x}\right),
\qquad x\in[-1,1].
\eq
In particular, $S_{[-1,1]}$ is a well defined
linear mapping acting from $C^{\iy}_{0}[-1,1]$ into $C^{\iy}[-1,1]$.

Let $\sigma:[-1,1]\to\R_{+}$ be a Lebesgue measurable, almost everywhere
nonzero function. We denote by $L^{p}_{\sigma}[-1,1]$
the space consisting of all Lebesgue mesurable functions
$f:[-1,1]\to\C$ for which
\bq
\|f\|_{L^{p}_{\sigma}[-1,1]}&:=&
\left(\int_{-1}^{1}\sigma^{p}(x)|f(x)|^{p}\,dx\right)^{1/p}\;\,<\;\,\iy.
\eq
For certain functions $\sigma$ it is possible to extend the
singular integral operator $S_{[-1,1]}$ as defined above on  
$C^{\iy}_{0}[-1,1]$
by continuity to a linear bounded operator acting on the
Banach space $L^{p}_{\sigma}[-1,1]$. 
The criteria is related to the $A_{p}$-condition on $[-1,1]$.

The following theorem was established in the case of the real line$\R$ rather than
the interval $[-1,1]$ first by Hunt, Muckenhoupt and Wheeden \cite{HMW}.
The theorem itself follows from the results of Coifman and Fefferman \cite{CF}.
For more information about the $A_{p}$ condition and the boundedness of the
singular integral operators on more general curves we refer \cite{BK}.

\begin{theorem}\label{t1.2}
Let $\sigma:[-1,1]\to\R_{+}$ be a Lebesgue measurable, almost everywhere
nonzero function. Assume that $\sigma\in L^{p}[-1,1]$ and $\sigma\iv\in  
L^{q}[-1,1]$,
where $1<p<\iy$, $1/p+1/q=1$.
Then $S_{[-1,1]}:C^{\iy}_{0}[-1,1]\to C^{\iy}[-1,1]$ can be continued  
by continuity to a linear bounded operator acting on  
$L^{p}_{\sigma}[-1,1]$ if and only if
$\sigma$ satisfies the $A_{p}$-condition on $[-1,1]$.
\end{theorem}
We remark in connection with the previous theorem that the assumptions
that $\sigma$ is nonzero almost everywhere and that $\sigma\in  
L^{p}[-1,1]$
imply that $C^{\iy}_{0}[-1,1]$ is  a dense linear subspace of  
$L^{p}_{\sigma}[-1,1]$.

The proof of this statement is similar to the proof of Lemma \ref{l2:1}  
below.
Let $\Ci$ stand for the set of all infinitely differentiable functions  
on $\T$,
and let $\Cio$ stand for the space of all $f\in\Ci$ such that
$f(t)$ and all of its derivatives vanish at $t=1$ and $t=-1$:
\bq
\Cio&=&\Big\{\;f\in\Ci\;:\;f^{(n)}(1)=f^{(n)}(-1)=0\mbox{ for all }
n\ge0\;\Big\}.
\eq

Let $\rho:\T\to\R_{+}$ be a Lebesgue measurable and almost everywhere
nonzero function.  We denote by $L^{p}_{\rho}(\T)$ the space of
all Lebesgue measurable functions $f:\T\to\C$ for which
\bq
\|f\|_{L^p_{\rho}(\T)} &:=&
\left(\frac{1}{2\pi}\int_{0}^{2\pi}
\rho^{p}(e^{i\theta})| f(e^{i\theta})|^{p}\,d\theta\right)^{1/p}
\;\,<\;\,\iy.
\eq
Let us remark that the dual space to $L^{p}_{\rho}(\T)$ can be  
identified
with $L^{q}_{\rho\iv}(\T)$ by means of the bilinear functional
\bq
\langle g,f\rangle &:=&\frac{1}{2\pi}
\int_{0}^{2\pi}\overline{g(e^{i\theta})}f(e^{i\theta})\,d\theta
\eq
with $f\in L^{p}_{\rho}(\T)$, $g\in L^{q}_{\rho\iv}(\T)$, $1<p<\iy$,  
$1/p+1/q=1$.

\begin{lemma}\label{l2:1}
Let $\rho\in\Lp$, $1<p<\iy$, and assume that $\rho$ is nonzero almost  
everywhere.
Then $\Cio$ is a dense subspace of  $L^{p}_{\rho}(\T)$.
\end{lemma}
\begin{proof}
We introduce the set
\bq
X&=& \Big\{\;\rho f\;:\;f\in\Cio\;\Big\}.
\eq
Obviously, $X\subseteq \Lp$, which implies that $\Cio$ is a subset of
$L^p_\rho(\T)$. The assertion that $\Cio$ is a dense subspace
in $L^{p}_{\rho}(\T)$ is
equivalent to the statement that $X$ is a dense subspace of $\Lp$.

We carry out the proof of this statement in several steps. First we  
prove that the closure of
$X$ contains all functions of the
form $f=\rho g$ with $g\in \Li$. Indeed, given $g\in \Li$ and
$\eps>0$, there are a subset $M\subset \T$ of Lebesgue measure less than
$\eps$ and a sequence $g_{n}\in\Cio$ such that $g_{n}\to g$ uniformly
on $\T\setminus M$ and $\|g_{n}\|_{\Li}\le \|g\|_{\Li}$. Now we can
estimate
\bq
\|\rho g-\rho g_{n}\|_{\Lp} &\le&
2\,\|\rho\|_{L^{p}(M)}\|g\|_{\Li}+
\|\rho\|_{\Lp}\|g_{n}-g\|_{L^{\iy}(\T\setminus M)}.\nn
\eq
Since $\|\rho\|_{L^{p}(M)}\to0$ as $\eps\to 0$ the assertion follows
easily.

Next we prove that $\Li$ is contained in the closure of $X$. Indeed,
given $f\in\Li$ we introduce the elements $g_{\eps}=\rho_{\eps}f$
where
\bq
\rho_{\eps}(t) &=& \left\{\ba{cl}
\rho\iv(t) & \mbox{ if } 0<\rho\iv(t)\le\eps\iv,\\[.5ex]
0 & \mbox{ if }\rho\iv(t)>\eps\iv.\ea\right.  \nn
\eq
Obviously, $\rho_{\eps}\in\Li$ and hence $g_{\eps}\in \Li$. Now we
estimate
\bq
\|\rho g_{\eps}-f\|_{\Lp} &=&
\|(\rho\rho_{\eps}-1)f\|_{\Lp}\;\,\le\;\,
\|\rho\rho_{\eps}-1\|_{\Lp}\|f\|_{\Li}.\nn
\eq
The function $1-\rho\rho_{\eps}$ is equal to the characteristic
function of
\bq
K_{\eps} &=&\Big\{\;t\in\T\;:\;\rho\iv(t)> \eps\iv\;\Big\}.\nn
\eq
Since $\rho$ is nonzero almost everywhere, the Lebesgue measure
of $K_{\eps}$ tends to zero as $\eps\to0$. Hence
$\|\rho\rho_{\eps}-1\|_{\Lp}\to 0$ as $\eps\to0$. It follows
that $\rho g_{\eps}$ approximates $f$.

Finally, it we note that $\Li$ is dense in $\Lp$.
\end{proof}

Let $Q=I-P$ stand for the complementary projection to the Riesz  
projection $P$.
Moreover, define the operators
\bq
G &=& \frac{1}{2}(I+J)(P-Q)(I-J),\\[.5ex]
G^* &=& \frac{1}{2}(I-J)(P-Q)(I+J).
\eq
We think of $G$ and $G^{*}$ as linear mappings acting from $\Cio$ into  
$\Ci$.
Notice in this connection that $P$ and $Q$ map $\Cio$ into $\Ci$.

\begin{proposition}\label{p2.2}
Let $1<p<\iy$ and $1/p+1/q=1$.
Assume that $\phi\in G\Li$ admits a weak asymmetric factorization
$\phi(t)=\phi_{-}(t)t^{\kappa}\phi_{0}(t)$ in $\Lp$.
Then the following is equivalent:
\begin{itemize}
\item[(i)]
The operator $B:=L(\phi_{0}\iv)(I+J)PL(\phi_{-}\iv):X_{1}\to X_{2}$
can be continued by continuity to a linear bounded operator acting
on $\Lp$.
\item[(iii)]
The operator $G^{*}:\Cio\to\Ci$ can be continued by continuity to a  
linear
bounded operator acting on $L^{q}_{\om\iv}(\T)$.
\item[(iii)]
The operator $G:\Cio\to\Ci$ can be continued by continuity to a linear
bounded operator acting on $L^{p}_{\omega}(\T)$.
\end{itemize}
Therein $\omega(t)=|\phi_{0}\iv(t)|$, and $X_{1}$ and $X_{2}$ are  
defined
by (\ref{f.X1}) and (\ref{f.X2}), respectively.
\end{proposition}
\begin{proof}
(i)$\Leftrightarrow$(ii):\
Assertion (i) is equivalent to the fact that there exists a constant
$M>0$ such that
\bq
\|L(\phi_{0}\iv)(I+J)PL(\phi_{-}\iv)f\|_{\Lp} &\le &
M\, \|f\|_{\Lp}\nn
\eq
for all $f\in X_{1}$. {}From the definition of $L^{p}_{\om}(\T)$ it
follows that the last inequality can be rewritten as
\bq
\|(I+J)PL(\phi_{-}\iv)f\|_{L^{p}_{\om}(\T)} &\le &
M\, \|f\|_{\Lp}.\nn
\eq
Since $\Cio$ is dense in $L^{q}_{\om\iv}(\T)$ by Lemma \ref{l2:1}, we  
obtain that
this is in turn equivalent to the statement that
\bq
|\langle g,(I+J)PL(\phi_{-}\iv)f\rangle| &\le& M\,
\|g\|_{L^{q}_{\om\iv}(\T)}\|f\|_{\Lp}\nn
\eq
for all $g\in\Cio$ and all $f\in X_{1}$. Next notice that
\bq
\langle g,(I+J)PL(\phi_{-}\iv)f\rangle  &=&
\langle P(I+J)g, L(\phi_{-}\iv)f\rangle\;\,=\;\,
\langle L((\phi_{-}^{-1})^{*})P(I+J)g,f\rangle\nn
\eq
by noting that $L(\phi_{-}\iv)f\in \Lq$. Hence the
above is equivalent to the statement that
\bq
|\langle L((\phi_{-}^{-1})^{*})P(I+J)g,f\rangle|
&\le& M\,\|g\|_{L^q_{\om\iv}(\T)}\|f\|_{\Lp}\nn
\eq
for all $g\in\Cio$ and all $f\in X_{1}$. Since  $X_1$ is dense in $\Lp$
we can reformulate this by saying that
\bq
\|L((\phi_-^{-1})^{*})P(I+J)g\|_{\Lq} &\le&
M\, \|g\|_{L^q_{\om\iv}(\T)}\nn
\eq
for all $g\in \Cio$. Because $\phi_-\iv(t)=\phi_0(t)\phi\iv(t)$
the latter can be rewritten as
\bq
\|P(I+J)g\|_{L^q_{\om\iv}(\T)} &\le&
M\, \|g\|_{L^q_{\om\iv}(\T)}.\nn
\eq
Since $\om(t)=\om(t\iv)$ the operator $(I+J)$ is bounded on
$L^q_{\om\iv}(\T)$. Moreover, since $P-Q=2P-I$
we can conclude that the latter is equivalent to
\bq
\|(P-Q)(I+J)g\|_{L^q_{\om\iv}(\T)}&\le&
M\, \|g\|_{L^q_{\om\iv}(\T)}
\eq
for all $g\in\Cio$. Noting that
\bq
G^{*}\;\,=\;\,
\frac{1}{2}(I-J)(P-Q)(I+J)&=&
(P-Q)(I+J)\nn
\eq
completes the proof of (i)$\Leftrightarrow$(ii).

(ii)$\Leftrightarrow$(iii):\
Since $\Cio$ is dense in both $L^p_{\om}(\T)$ and
$L^q_{\om\iv}(\T)$ by Lemma \ref{l2:1}, it is easily seen that (ii) is  
equivalent to
the statement that
\bq
|\langle G^{*}g,f\rangle|&\le& M\,
\|g\|_{L^p_{\om}(\T)}\|f\|_{L^q_{\om\iv}(\T)}\nn
\eq
for all $f,g\in\Cio$. Moreover, (iii) is equivalent to
the statement that
\bq
|\langle g,Gf\rangle|&\le& M\,
\|g\|_{L^p_{\om}(\T)}\|f\|_{L^q_{\om\iv}(\T)}\nn
\eq
for all $f,g\in\Cio$. Since
$\langle G^{*}g,f\rangle=\langle g,Gf\rangle$ the result follows.
\end{proof}

Let $C^{\iy}_{\even}(\T)$ stand for the space of all functions $f\in\Ci$
which are even, i.e., for which $f(t)=f(t\iv)$, $t\in\T$. Moreover,  
introduce the operators
\bq
U &:&
\hat{f}(x)\in C^\iy_{0}[-1,1]
\mapsto f(e^{i\theta}):=\hat{f}(\cos\theta)\in\Cio,\\[1ex]
V &:&
f(t)\in C^{\iy}_{\even}(\T)\mapsto
\hat{f}(\cos\theta):=f(e^{i\theta})\in C^{\iy}[-1,1],
\eq
and let
\bq
\chi(e^{i\theta}) &=&
\left\{\ba{cl} 1 &\mbox{ if }0<\theta<\pi,\\[.5ex]
  -1 &\mbox{ if }\pi<\theta<2\pi.\ea\right.
\eq
Clearly, the image of $U$ is the set of even functions defined on $\T$.

\begin{proposition}\label{p2.3}
For $\hat{f}\in C^\iy_{0}[-1,1]$ we have
\bq
S_{[-1,1]}\hat{f} &=& \frac{1}{2}\,V  L((1+t\iv)\iv)GL(\chi (1+t\iv)) U  
\hat{f}.
\eq
\end{proposition}
\begin{proof}
Given $\hat{f}\in C^\iy_0[-1,1]$, we introduce the functions
$$
f=L(\chi(1+t\iv)) U \hat{f},\quad   g=\frac{1}{2}\,Gf,\quad
\hat{g}=V  L((1+t)\iv) g.
$$
Notice that
$$
G=\frac{1}{2}(I+J)(P-Q)(I-J)=(P-Q)(I-J).
$$
Hence $g(t)=t\iv g(1/t)$ and it follows that $L((1+t\iv)\iv)g$
is an even function. Moreover, it is easily seen that
$Jf= -f$ whence it follows that $g=\frac{1}{2}\,Gf=(P-Q)f$. It is well  
known
that the singular integral operator $S=P-Q$ on $\T$ can be written as
\bq
(Sf)(t) &=& \frac{1}{\pi i}\int_\T \frac{f(s)-f(t)}{s-t}\,ds+f(t)\nn
\eq
for functions $f\in\Ci$. From this we deduce the relations
\bq
f(e^{i\theta}) &=& (1+e^{-i\theta})\chi(e^{i\theta})\hat{f}(\cos\theta),
\nn\\[1ex]
g(e^{i\theta}) &=&
\frac{1}{\pi}\int_{-\pi}^{\pi}  
\frac{f(e^{i\varphi})-f(e^{i\theta})}{1-e^{i(\theta-\varphi)}}\,d\varphi
+f(e^{i\theta}),\nn\\[1ex]
\hat{g}(\cos\theta) &=& (1+e^{-i\theta})\iv g(e^{i\theta}).\nn
\eq

We split the integral appearing in the second equation
into two parts integrating on $[0,\pi]$ and $[-\pi,0]$, respectively,
and make a change of variables $\varphi\mapsto -\varphi$ in the second  
integral.
This gives
\bq
g(e^{i\theta}) &=&
\frac{1}{\pi}\int_{0}^{\pi} \left(  
\frac{f(e^{i\varphi})-f(e^{i\theta})}{1-e^{i(\theta-\varphi)}}
+\frac{f(e^{-i\varphi})-f(e^{i\theta})}{1- 
e^{i(\theta+\varphi)}}\right)\,d\varphi
+f(e^{i\theta}).\nn
\eq
Since $f(e^{-i\varphi})=-e^{i\varphi}f(e^{i\varphi})$ and since
\bq
\frac{1}{1-e^{i(\theta-\varphi)}}
-\frac{e^{i\varphi}}{1-e^{i(\theta+\varphi)}}
&=&
\frac{(1+e^{-i\theta})(e^{i\varphi}-1)}{2(\cos\varphi-\cos\theta)},
\nn\\
\frac{1}{1-e^{i(\theta-\varphi)}}
+\frac{1}{1-e^{i(\theta+\varphi)}}
&=&  
\frac{e^{i\varphi}+e^{-i\varphi}-2e^{-i\theta}}{2(\cos\varphi- 
\cos\theta)},
\nn
\eq
it follows that
\bq
g(e^{i\theta}) &=&
\frac{1}{\pi}\int_{0}^{\pi} \left(
\frac{(1+e^{-i\theta})(e^{i\varphi}-1)f(e^{i\varphi})}{2(\cos\varphi- 
\cos\theta)}
-\frac{(e^{i\varphi}-e^{-i\varphi})f(e^{i\theta})}{2(\cos\varphi- 
\cos\theta)}\right)\,d\varphi
\nn\\[1ex]
&& -\frac{f(e^{i\theta})}{\pi}\int_{0}^{\pi}
\frac{e^{-i\varphi}-e^{-i\theta}}{\cos\varphi-\cos\theta}\,d\varphi
+f(e^{i\theta}).\nn
\eq

If we assume $0<\theta<\pi$, we obtain
\bq
g(e^{i\theta}) &=&
\frac{1}{\pi}\int_{0}^{\pi}  
\frac{(1+e^{-i\theta})(e^{i\varphi}-e^{-i\varphi})(\hat{f}(\cos\varphi)-
\hat{f}(\cos\theta))}{2(\cos\varphi-\cos\theta)}
\,d\varphi
\nn\\[1ex]
&& -\frac{f(e^{i\theta})}{\pi i}\int_{0}^{\pi}
\frac{\sin\varphi-\sin\theta}{\cos\varphi-\cos\theta}\,d\varphi.\nn
\eq
The first integral is equal to $(1+e^{-i\theta})$ times
\bq
\frac{i}{\pi}\int_{0}^{\pi} \frac{(\hat{f}(\cos\varphi)-
\hat{f}(\cos\theta))\sin\varphi}{\cos\varphi-\cos\theta}
\,d\varphi
&=&
\frac{1}{\pi i}\int_{-1}^1  
\frac{\hat{f}(y)-\hat{f}(\cos\theta)}{y-\cos\theta}\,dy.\nn
\eq
The second integral equals $f(e^{i\theta})$ times
$$
\frac{1}{\pi i}\int_{0}^{\pi}
\cot\left(\frac{\varphi+\theta}{2}\right)\,d\varphi
=
\left[\frac{2}{\pi  
i}\ln\sin\left(\frac{\varphi+\theta}{2}\right)\right]_{\varphi=0}^{\pi}
=
\frac{1}{\pi i}\ln \left(\frac{1-\cos\theta}{1+\cos\theta}\right).
$$
Putting the pieces together we obtain
\bq
g(e^{i\theta}) &=&
\frac{(1+e^{-i\theta})}{\pi i}\int_{-1}^1
\frac{\hat{f}(y)-\hat{f}(\cos\theta)}{y-\cos\theta}\,dy
+
\frac{f(e^{i\theta})}{\pi i}
\ln \left(\frac{1-\cos\theta}{1+\cos\theta}\right).\nn
\eq

Note that the above assumption $0<\theta<\pi$ is not an essential
restriction since $\hat{g}$ is determined by the formula
\bq
\hat{g}(\cos\theta) &=& (1+e^{-i\theta})\iv g(e^{i\theta}),\qquad  
0<\theta<\pi.\nn
\eq
It follows that
\bq
\hat{g}(x) &=& \frac{1}{\pi  
i}\int_{-1}^1\frac{\hat{f}(y)-\hat{f}(x)}{y-x}\,dy+
\frac{\hat{f}(x)}{\pi i}\ln\left(\frac{1-x}{1+x}\right).\nn
\eq
This is equal to the singular integral operator $S_{[-1,1]}$ applied
to the function $\hat{f}$. Hence $\hat{g}=S_{[-1,1]}\hat{f}$ which is
the assertion.
\end{proof}

Now we are able to present the proof of Theorem \ref{t1.3}.

\begin{proofof}{Theorem \ref{t1.3}}
It is obvious from Proposition \ref{p1.1} and Theorem \ref{t1.1}
that the Fredholmness of $M(\phi)$ implies assertions (a) and (b).

Hence it is sufficient to prove the following. If conditions
(a) and (b) are fulfilled, then $M(\phi)$ is Fredholm if and only if
$\sigma$ satisfies the $A_p$-condition.

Under these assumptions we deduce from Theorem \ref{t1.1}
that the Fredholmness of $M(\phi)$ on $\Hp$ is equivalent to
the existence of a bounded continuation of the operator
$B=L(\phi_0\iv)(I+J)PL(\phi_-\iv):X_1\to X_2$ on $\Lp$.
We apply Proposition \ref{p2.2} and see that this existence
is equivalent to the condition that
\bq
\|Gg\|_{L^p_{\om}(\T)}&\le&M\,
\|g\|_{L^p_{\om}(\T)}\nn
\eq
for all $g\in \Cio$ where $\om(t):=|\phi_0\iv(t)|$. 

Next we claim that this, in turn, is equivalent to the condition
that
\bq
\|Gg\|_{L^p_{\om}(\T)}&\le&M\,
\|g\|_{L^p_{\om}(\T)}\nn
\eq
for all $g\in \Cio$ for which $Jg=-g$. In order to prove the non-trivial
part of this equivalence, we decompose an arbitraryly given
$g\in\Cio$ into $g=g_1+g_2$ where $g_1=\frac{1}{2}(I+J)g$ and
$g_2=\frac{1}{2}(I-J)g$. The function $g_1$ lies in the kernel of $G$
while $g_2(t)=(g(t)-t\iv g(t\iv))/2$ belongs to $\Cio$ and
satisfies the relation $Jg_2=-g_2$. Moreover, since $\om$ is an
even function, the operator $(I-J)/2$ is bounded on
$L^p_{\om}(\T)$.
We obtain the estimate
$$
\|Gg\|_{L^{p}_{\om}(\T)} = 
\|Gg_{2}\|_{L^{p}_{\om}(\T)}\le
M\,\|g_{2}\|_{L^{p}_{\om}(\T)}\le
M\,\|g\|_{L^{p}_{\om}(\T)},
$$
proves the claim.

Next we remark that the operator $L(\chi(1+t\iv))U$
maps the space $C^\iy_0[-1,1]$ onto the subspace of functions
$g\in \Cio$ satisfying $Jg=-g$. This allows us to make the
substitution $g=L(\chi(1+t\iv))Uf$ with $f\in C^{\iy}_{0}[-1,1]$.
We obtain the equivalent estimate
\bq
\|GL(\chi(1+t\iv))Uf\|_{L^p_{\om}(\T)}&\le&M\,
\|L(\chi(1+t\iv))Uf\|_{L^p_{\om}(\T)}\nn
\eq
for all $f\in C^\iy_0[-1,1]$. Clearly, the last estimate can be written
in the form
\bq
\|V  L((1+t\iv)\iv)GL(\chi(1+t\iv)Uf\|_{L^p_{\sigma}[-1,1]}
&\le & M\,\|f\|_{L^p_{\sigma}[-1,1]}\nn
\eq
for all $f\in C^\iy_0[-1,1]$, where
\bq
\sigma(\cos\theta)&=&
\frac{\omega(e^{i\theta})|1+e^{-i\theta}|}{\sqrt{2}\,|\sin\theta|^{1/p}}
\;\,=\;\,
|\phi_0\iv(e^{i\theta})|\,
\frac{|1+\cos\theta|^{1/2q}}{|1-\cos\theta|^{1/2p}}.\nn
\eq
Along with Proposition \ref{p2.3} and Theorem \ref{t1.2} this
completes the proof.
\end{proofof}


\section{Application to piecewise continuous functions}

We now apply the previous results in order to obtain
necessary and sufficient conditions for the operator $M(\phi)$
to be invertible or Fredholm on $\Hp$ for piecewise continuous
functions $\phi$ with finitely many jumps. These results have already  
been
established in \cite{BaEh2} (and in \cite{BaEh1} for the case $p=2$).
The proofs given in \cite{BaEh1,BaEh2} relies on
the results establish in \cite{RS} by help of Banach algebra methods.
The proof which we will give here is more direct and relies entirely
on the factorization methods developed here and in \cite{BaEh2}
in connection with the $A_p$-condition.

We restrict to piecewise continuous functions with a
{\em finite} number of discontinuities because these
functions can be written in a convenient manner which is useful in
many instances. It is well known that any piecewise continuous
and nonvanishing function with a finite number of discontinuities
at the points $\theta_{1},\dots, \theta_{R}$
can be written as a product
\bq\label{f3.23}
\phi(e^{i\theta}) &=&
b(e^{i\theta})\prod_{r=1}^{R}t_{\beta_{r}}(e^{i(\theta-\theta_{r})})
\eq
where $b$ is a nonvanishing continuous function and
\bq\label{f3.24}
t_{\beta}(e^{i\theta}) &=&
\exp(i\beta(\theta - \pi)),\qquad 0 < \theta < 2\pi.
\eq
Notice that the parameters $\beta_{1},\dots,\beta_{R}$ in this formula  
are
uniquely determined up to an additive integer. In fact,
\bq\label{tbeta1}
\frac{\phi(e^{i\theta_{r}}-0)}{\phi(e^{i\theta_{r}}+0)}
&=&
\frac{t_{\beta_{r}}(1-0)}{t_{\beta_{r}}(1+0)}
\;\,=\;\,
\exp(2\pi i\beta_{r}).
\eq
Moreover, the formula
\bq\label{tbeta2}
t_{\beta+n}(t) &=&
(-t)^{n}t_{\beta}(t),\quad t\in\T,
\eq
holds for $n\in\Z$.

The parameters in the representation (\ref{f3.23}) are useful to decide  
Fredholmness and
invertibility. For example, the Toeplitz operators $T(\phi)$ with a  
piecewise continuous
function $\phi$ with finitely many jump discontinuities
is invertible on $\Hp$
if and only if the function $\phi$ can be represented in the form  
(\ref{f3.23}) with
$-1/q<\Re \beta _{r} < 1/p$ and  the winding number of $b$ equal to  
zero.
If the operator is Fredholm with Fredholm index $\kappa$,
then the zero is replaced by $-\kappa$.
(\cite{BS1}, Chapter 6)

Before stating the analogue of this result for Toeplitz plus Hankel  
operators $M(\phi)$
we have to establish the following auxiliary result.

\begin{lemma}\label{l3.1}
Let $1<p<\iy$, $1/p+1/q=1$, $-1=x_{0}<x_{1}<\dots<x_{R}<x_{R+1}=1$ and
\bq
\sigma(x) &=&\prod_{r=0}^{R+1}|x-x_{r}|^{\alpha_{r}}
\eq
If $-1/p<\alpha_{r}<1/q$ for each $0\le r\le R+1$, then $\sigma\in  
L^{p}[-1,1]$,
$\sigma\iv\in L^{q}[-1,1]$ and $\sigma$ satisfies the $A_{p}$-condition.
\end{lemma}
\begin{proof}
This can be verified straightforwardly.
\end{proof}

The promised result is the following:

\begin{theorem} \label{t6.1}
Let $1<p<\iy$, $1/p+1/q=1$.
Suppose that $\phi$ has finitely many jump discontinuities. Then
$M(\phi)$ on is Fredholm on $\Hp$
if and only if $\phi$ can be written in the form
\bq\label{f3.25}
\phi(e^{i\theta}) &=&
  b(e^{i\theta})  t_{\beta^{+}}(e^{i\theta})
  t_{\beta^{-}}(e^{i(\theta-\pi)}) \prod_{r=1}^{R}
  t_{\beta_r^+}(e^{i(\theta -\theta_{r})})
  t_{\beta_r^-}(e^{i(\theta +\theta_{r})})
\eq
where $b$ is a continuous nonvanishing function on $\T$,
the numbers $\theta_{1}, \dots, \theta_{R} \in (0, \pi)$
are distinct, and
\begin{enumerate}
\item[(i)] $-1/q < \Re(\beta_r^+ + \beta_r^-) < 1/p$
for each $1\le r\le R$,
\item[(ii)] $ -1/2 -1/2q < \Re\beta^{+} < 1/2p$ and $
-1/2q < \Re\beta^{-} < 1/2 +1/2p .$
\end{enumerate}
Moreover, in this case,
\bq\label{f3.ker-coker}
\dim\ker M(\phi)=\max\{0,-\wind(b)\},\qquad
\dim\ker M^*(\phi)=\max\{0,\wind(b)\}.
\eq
\end{theorem}
\begin{proof}
In the first step we prove that $M(\psi)$ is a Fredholm operator on
$\Hp$ with Fredholm index zero if  $\psi$ is of the form
\bq
\psi(e^{i\theta}) &=&
  t_{\beta^{+}}(e^{i\theta})
  t_{\beta^{-}}(e^{i(\theta-\pi)}) \prod_{r=1}^{R}
  t_{\beta_r^+}(e^{i(\theta -\theta_{r})})
  t_{\beta_r^-}(e^{i(\theta +\theta_{r})})
\eq
and the parameters satisfy the conditions (i) and (ii). In regard to  
Theorem
\ref{t1.3} it suffices to construct a weak asymmetric factorization
of $\psi$ and to prove that the corresponding weight $\sigma$ satisfies
the $A_p$-condition. For this purpose we introduce the
functions
$$
\eta_\beta(t)=(1-t)^\beta,\qquad
\xi_\beta(t)=(1-t\iv)^\beta.
$$
Notice that $t_\beta(t)=\eta_\beta(t)\xi_{-\beta}(t)$.
Then we can factor $\psi(t)=\psi_-(t)\psi_0(t)$ with
\bq
\psi_-(e^{i\theta}) &=& \xi_{-2\beta^+}(e^{i\theta})
\xi_{-2\beta^-}(e^{i(\theta-\pi)})\nn\\[1ex]
&&\times \prod_{r=1}^R
\xi_{-\beta_r^+-\beta_r^-}(e^{i(\theta-\theta_r)})
\xi_{-\beta_r^+-\beta_r^-}(e^{i(\theta+\theta_r)}),
\nn\\[1ex]
\psi_0(e^{i\theta}) &=&
\eta_{\beta^+}(e^{i\theta}) \xi_{\beta^+}(e^{i\theta})
\eta_{\beta^-}(e^{i(\theta-\pi)}) \xi_{\beta^-}(e^{i(\theta-\pi)})
\nn\\[1ex]
&&\times \prod_{r=1}^R
\eta_{\beta_r^+}(e^{i(\theta-\theta_r)})
\xi_{\beta_r^+}(e^{i(\theta+\theta_r)})
\eta_{\beta_r^-}(e^{i(\theta+\theta_r)})
\xi_{\beta_r^-}(e^{i(\theta-\theta_r)}).\nn
\eq
Because of conditions (i) and (ii), it can be checked straightforwardly  
that the function
\bq
(1+e^{-i\theta})\psi_-(e^{i\theta}) &=&
\xi_{-2\beta^+}(e^{i\theta})
\xi_{-2\beta^-+1}(e^{i(\theta-\pi)})\nn\\[1ex]
&&\times \prod_{r=1}^R
\xi_{-\beta_r^+-\beta_r^-}(e^{i(\theta-\theta_r)})
\xi_{-\beta_r^+-\beta_r^-}(e^{i(\theta+\theta_r)})\nn
\eq
belongs to $\ovl{\Hp}$ and the function
\bq
(1-e^{-i\theta})\psi_-\iv(e^{i\theta}) &=&
\xi_{2\beta^++1}(e^{i\theta})
\xi_{2\beta^-}(e^{i(\theta-\pi)})\nn\\[1ex]
&&\times \prod_{r=1}^R
\xi_{\beta_r^++\beta_r^-}(e^{i(\theta-\theta_r)})
\xi_{\beta_r^++\beta_r^-}(e^{i(\theta+\theta_r)})\nn
\eq
belongs to $\ovl{\Hq}$.
{}From the fact that $\psi_0$ is even and that
$\psi_0(t)=\psi_-(t)\iv\psi(t)$, it follows that
the function $\psi_0$ fulfills all the
necessary conditions in regard to a weak asymmetric factorization.
Hence $\psi(t)=\psi_-(t)\psi_0(t)$ is indeed a weak asymmetric
factorization with index zero.

In order to calculate the corresponding weight function
(\ref{f1.sigma})  consider
\bq
\psi_{0}(e^{i\theta}) &=&
|1-e^{i\theta}|^{2\beta^{+}}|1+e^{i\theta}|^{2\beta^{-}}
\nn\\[1ex]
&&\times
\prod_{r=1}^{R}
|1-e^{i(\theta-\theta_{r})}|^{\beta_{r}^{+}+\beta_{r}^{-}}
|1-e^{i(\theta+\theta_{r})}|^{\beta_{r}^{+}+\beta_{r}^{-}}
t_{\frac{\beta_{r}^{+}-\beta_{r}^{-}}{2}}(e^{i(\theta-\theta_{r})})
t_{\frac{\beta_{r}^{-}-\beta_{r}^{+}}{2}}(e^{i(\theta+\theta_{r})}),\nn
\eq
and observe that  
$|1-e^{i\theta}|=(2-2\cos\theta)^{1/2}=2|\sin(\frac{\theta}{2})|$
and  
$2\sin\frac{\theta- 
\theta_{r}}{2}\sin\frac{\theta+\theta_{r}}{2}=\cos\theta_{r}- 
\cos\theta$.
Hence
\bq
\psi_{0}\iv(e^{i\theta}) &=&
\sigma_0(\cos\theta)(1-\cos\theta)^{-\beta^+}(1+\cos\theta)^{-\beta^-}
\prod_{r=1}^R
|\cos\theta-\cos\theta_r|^{-\beta_r^+-\beta_r^-},\nn
\eq
where $\sigma_{0}(x)\in G\Li$ is a functions which comes from  
collecting the terms
$t_{\frac{\beta_{r}^{+}-\beta_{r}^{-}}{2}}(e^{i(\theta-\theta_{r})})$,
$t_{\frac{\beta_{r}^{-}-\beta_{r}^{+}}{2}}(e^{i(\theta+\theta_{r})})$  
and certain constants.
It follows that $\sigma$ evaluates to
\bq
\sigma(x) &=&
|\sigma_0(x)|(1-x)^{-\Re\beta^+-1/2p}(1+x)^{-\Re\beta^-+1/2q}
\prod_{r=1}^R
|x-\cos\theta_r|^{-\Re\beta_r^+-\Re\beta_r^-}.
\eq
It suffices to apply Lemma \ref{l3.1} in order to see that $\sigma$  
satisfies the
$A_{p}$-condition.

In the second step we prove that $M(\phi)$ is a Fredholm operator if  
$\phi$ is given by
(\ref{f3.25}) with conditions (i) and (ii) being fulfilled and if the  
function $b$ is continuous and nonvanishing. We can write
$\phi(t)=b(t)\psi(t)$ where $\psi$ is as above. {}From well-known  
identities
for Toeplitz and Hankel operators,
\bq
T(\phi) &=& T(b)T(\psi)+H(b)H(\wt{\psi}),\nn\\[1ex]
H(\phi) &=& T(b)H(\psi)+H(b)T(\wt{\psi}),\nn
\eq
where $\wt{\psi}(t)=\psi(t\iv)$, it follows that
\bq
M(\phi) &=& T(b)M(\psi)+H(b)M(\wt{\psi})\nn
\eq
Under the assumption on $b$ the operator $H(b)$ is compact and the  
operator
$T(b)$ is Fredholm with Fredholm index equal to $-\wind(b)$. Since we  
have just proved
that $M(\psi)$ is Fredholm with Fredholm index zero, it follows
that $M(\phi)$ is Fredholm with Fredholm  index equal to
$-\wind(b)$.

Hence we have proved the ``if'' part of the theorem and also computed  
the
Fredholm index of $M(\phi)$. Now we apply Theorem \ref{t1.1} with  
formula
(\ref{ker-coker}). This formula implies that the Fredholm index is  
equal to $-\kappa$,
where $\kappa$ is the index of the asymmetric factorization of $\phi$.  
Hence
$\kappa=\wind(b)$ and formula (\ref{f3.ker-coker}) follows.

In the last step we are going to prove the ``only if'' part of the  
theorem.
It is settled by a well-known perturbation argument.
Suppose that $M(\phi)$ is a Fredholm operator with index $\kappa$, say.
We conclude from Proposition \ref{p1.1} that $\phi\in G\Li$. Since  
$\phi$ has only a
finite number of jump discontinuities, this implies that
$\phi$ can be written in the form (\ref{f3.23}) or in the form
(\ref{f3.25}) if we put some of the $\beta$-parameters equal to zero
if necessary. Therein the function $b$ is continuous and nonvanishing
on $\T$. Moreover, due to formula (\ref{tbeta1}) and (\ref{tbeta2})
we can choose the $\beta$-parameters to satisfy the conditions
\begin{enumerate}
\item[(i*)] $-1/q < \Re(\beta_r^+ + \beta_r^-) \le 1/p$
for each $1\le r\le R$,
\item[(ii*)] $ -1/2 -1/2q < \Re\beta^{+} \le 1/2p$ and $
-1/2q < \Re\beta^{-} \le 1/2 +1/2p .$
\end{enumerate}
Assume contrary to what we want to prove, namely,
that in at least one instance we have equality in
the above conditions (i*) and (ii*). We are going to perturbate the  
$\beta$-parameters (and thus the
function $\phi$) in two different ways in order to arrive at  a  
contradiction. Remark that
since $M(\phi)$ is assumed to be Fredholm, the Fredholm index is  
constant with respect to
any small perturbation. We first perturbate by replacing all  
$\beta$-parameters
by $\beta-\eps$ where $\eps>0$ is sufficiently small. This turns the  
conditions (i*) and (ii*) into
(i) and (ii). Applying the ``if'' part of the theorem with the formula  
for the Fredholm index, it follows
that $\kappa=-\wind(b)$. In the second perturbation we do the same  
substitution
except for one of the instances of equality in (i*) and (ii*) where we  
replace the corresponding
$\beta^{\pm}$ by $\beta^{\pm}-1+\eps$, or, the corresponding
$\beta_{r}^{+}+\beta_{r}^{-}$ by $\beta_{r}^{+}+\beta_{r}^{-}-1+\eps$,  
respectively.
Moreover, we have to replace $b(t)$ by $t\, b(t)$ times a certain  
constant due to formula
(\ref{tbeta2}). This is again a small perturbation of $\phi$, which  
leaves the Fredholm index
unchanged. The corresponding parameters fulfill (i) and (ii), but we  
obtain
$\kappa=-\wind(t\, b(t))=-1-\wind(b)$ contradicting the above formula.
This completes the proof of the ``only if'' part.
\end{proof}


\end{document}